\documentclass[draft,12pt]{amsart}
\usepackage{amssymb}

\theoremstyle{plain}
\newtheorem{lemma}{Lemma}[section]
\newtheorem{proposition}[lemma]{Proposition}
\newtheorem{theorem}[lemma]{Theorem}
\newtheorem{corollary}[lemma]{Corollary}
\newtheorem{conjecture}[lemma]{Conjecture}

\newtheorem{problem}[lemma]{Problem}

\theoremstyle{definition}

\newtheorem{example}[lemma]{Example}

\newtheorem{remark}[lemma]{Remark}

\newcommand{\ggcaffil}[1]{\dedicatory{\textup{\larger{#1}}}}

\newcommand{\ggcqedsymbol}{$\square$}
\newcommand{\ggcqed}{\hbox{}\nobreak\hbox{\quad\ggcqedsymbol}}
\newcommand{\ggcendpf}{\ggcqed}
\newcommand{\ggcnopf}{\ggcqed}

\newcommand{\ggcenddef}{\ggcqed}
\newcommand{\ggcendeg}{\ggcenddef}
\newcommand{\ggcendrem}{\ggcenddef}

\newcommand{\bZ}{\ensuremath{\mathbb{Z}}}
\newcommand{\iv}[2]{\ensuremath{\left[#1,#2\right]}}
\newcommand{\gz}[1]{\ensuremath{\bZ(#1)}}
\newcommand{\gzmatch}[1]{\ensuremath{\bZ\left(#1\right)}}
\newcommand{\gzd}{\gz{D}}
\newcommand{\chd}{\ensuremath{\chi(D)}}
\newcommand{\gcdd}{\ensuremath{\gcd(D)}}

\begin{document}
\title{Coloring Distance Graphs on the Integers}
\author{Glenn G. Chappell}
\ggcaffil{Department of Mathematics, Southeast Missouri
  State University}
\address{Department of Mathematics\\
  Southeast Missouri State University\\
  Cape Girardeau, MO 63701\\
  USA}
\email{gchappell@semovm.semo.edu}
\date{April 23, 1998}
\begin{abstract}
  Given a set $D$ of positive integers, the associated
  \emph{distance graph} on the integers is the graph with the integers as
  vertices and an edge between distinct vertices if their difference lies
  in $D$.
  We investigate the chromatic numbers of distance graphs.
  We show that, if $D=\left\{ d_1,d_2,d_3,\dotsc\right\}$,
  with $d_n\mid d_{n+1}$ for all $n$, then the distance graph has a
  proper 4-coloring. We further find the exact chromatic numbers
  of all such distance graphs.
  Next, we characterize those distance graphs that have periodic
  proper colorings and show a relationship between the chromatic
  number and the existence of periodic proper colorings.
\end{abstract}
\maketitle

\section{Introduction} \label{S:intro}

What is the least number of classes into which the integers
can be partitioned, so that no two members of the same class
differ by a square? What if ``square'' is replaced by ``factorial''?

Questions like these can be formulated as graph coloring problems.
Given a set $D$ of positive integers,
the \emph{distance graph} $\gzd$ is the graph with the integers
as vertices and an edge between distinct vertices if
their difference lies in $D$;
we call $D$ the \emph{distance set} of this graph.
A \emph{proper coloring} of a graph is an assignment of colors
to the vertices so that no two vertices joined by an edge
receive the same color.
The \emph{chromatic number} of a graph $G$, denoted by $\chi(G)$,
is the least number of colors in a proper coloring.
We abbreviate $\chi\big(\gzd\big)$ by $\chd$.
We refer to~\cite{BoMu76,WesD96} for
graph-theoretic terminology not defined here.

When $D$ is the set of all positive squares, we call $\gzd$
the \emph{square distance graph}.
When $D$ is the set of all factorials, we obtain the
\emph{factorial distance graph}.
The questions at the beginning of this section ask for
the chromatic numbers of these two graphs.
We will study the chromatic numbers of these and other
distance graphs on the integers.

Distance graphs on the integers were introduced by
Eggleton, Erd{\H o}s, and Skilton in~\cite{EES85}.
In \cite{EES85,EES86}, the problem was posed of
characterizing those distance sets $D$,
containing only primes, such that $\chd=4$.
This problem was studied in
\cite{EggR88,EES90,VoWa91,VoWa94};
see also \cite{EES88}.
More recently,
\cite{CCH97,DeZh97} have discussed
the chromatic numbers of more general distance graphs with
distance sets having 3 or 4 elements.

In this paper, we are primarily interested in distance
graphs for which the distance set is infinite,
although our results apply to finite distance sets as well.
We begin in Section~\ref{S:easy}
with some easy lemmas on connectedness and bounds on
the chromatic number.
In Section~\ref{S:dc}, we consider distance graphs for which the
distance set is totally ordered by the divisibility relation.
We determine the chromatic numbers of all such graphs; in particular,
we prove that they are all $4$-colorable.
In Section~\ref{S:pc}, we study periodic proper colorings of
distance graphs and their relationship to the chromatic number.

Throughout this paper
we will use standard notation for intervals to denote
sets of consecutive integers.
For example, $\iv{2}{6}$ denotes the set
$\left\{2,3,4,5,6\right\}$.

\section{Basic Results} \label{S:easy}

In this section, we establish some basic facts about the
connectedness and chromatic number of distance graphs.
The results of this section have all been at least partially
stated in earlier works.

Our first result characterizes those
distance sets for which the distance graph is connected.
This result has been partially stated or implicitly
assumed in a number of earlier works;
see~\cite[p.~95]{EES85}.
For $D$ a set of positive integers,
we note that $\gcdd$ is well defined when $D$ is
infinite.
Given a real number $k$ and a set $D$, we denote by
$k\cdot D$ the set $\left\{\,kd:d\in D\,\right\}$.

\begin{lemma} \label{L:conn}
Let $D$ be a nonempty set of positive integers.
The graph $\gzd$ is connected if
and only if $\gcdd=1$.
Further, each component of $\gzd$ is isomorphic to
$\gzmatch{\frac{1}{\gcdd}\cdot D}$.
\end{lemma}

\begin{proof}
There is a path between vertices $k$ and $k+1$ if and only if
there exist $d_1,\dotsc,d_a,e_1,\dotsc,e_b\in D$ such that
$d_1+\dotsb+d_a-e_1-\dotsb-e_b=1$. This happens
precisely when $\gcdd=1$, and so the first
statement of the lemma is true.

For the second statement, one isomorphism is the function\break
$\varphi\colon\gcdd\cdot\bZ\to\bZ$ defined by
$\varphi(k)=\frac{k}{\gcdd}$.\ggcendpf
\end{proof}

When we determine the chromatic numbers of distance graphs,\break
Lemma~\ref{L:conn} will often allow us to assume that the GCD of
the distance set is $1$.

Next, we prove a useful upper bound on the chromatic number.
This result is a slight generalization of a result of
Chen, Chang, and Huang \cite[Lemma~2]{CCH97}.

\begin{lemma} \label{L:nomult}
Let $D$ be a nonempty set of positive integers, and
let $k$ be a positive integer.
If $\frac{1}{\gcdd}\cdot D$ contains no multiple of $k$,
then $\chd\le k$.
\end{lemma}

\begin{proof}
Let $D$ and $k$ be as stated.
By Lemma~\ref{L:conn} we may assume that $\gcdd=1$.
Thus, we assume that $D$ contains no multiple of $k$.
We color the integers with colors $\iv{0}{k-1}$,
assigning to each integer $i$ the color corresponding
to the residue class of $i$ modulo $k$.
Two integers will be assigned the same color precisely when
they differ by a multiple of $k$.
Since no multiple of $k$ occurs in $D$, this is a proper
$k$-coloring of $\gzd$.\ggcendpf\end{proof}

The converse of Lemma~\ref{L:nomult} holds when $k=2$.
This gives us a characterization of bipartite
distance graphs:
$\gzd$ is bipartite precisely when
$\frac{1}{\gcdd}\cdot D$ contains no multiple of $2$,
that is, when all elements of $D$ have the same power
of $2$ in their prime factorizations.
This result has been partially stated in earlier works;
see~\cite[Thms.~8 \&~10]{EES85} and~\cite[Thms.~3 \&~4]{CCH97}.

\begin{proposition} \label{P:chi2}
Let $D$ be a set of positive integers.
The graph $\gzd$ is bipartite
if and only if there exists a
non-negative integer $k$ so that
$\frac{1}{2^k}\cdot D$ contains only odd integers.
\end{proposition}

\begin{proof}
We may assume $D\ne\emptyset$.
Since a graph is bipartite if and only if each component is
bipartite, we may also assume,
by Lemma~\ref{L:conn}, that $\gcdd=1$.
For such $D$ we show that $\gzd$ is bipartite
if and only if each element of $D$ is odd.

$(\Longrightarrow)$ Since $\gcdd=1$, $D$ must have
an odd element $d$.
Suppose that $D$ has an even element $e$.
If we begin at $0$, take $e$ steps
in the positive direction,
each of length $d$, ending at $de$, and then take $d$ steps
in the negative direction,
each of length $e$, ending at $0$, then we have followed a closed
walk of odd length.
Formally, the set
\[\left\{0,d,d\cdot2,\dotsc,d(e-1),de,(d-1)e,(d-2)e,\dotsc,2e,e\right\}\]
is the vertex set of an odd circuit, and so $\gzd$ is
not bipartite.

$(\Longleftarrow)$
If every element of $D$ is odd, then $\chd\le2$,
by Lemma~\ref{L:nomult}.\ggcendpf\end{proof}

The converse of Lemma~\ref{L:nomult} does not hold when
$k>2$.
For example, let $k>2$, and let $D=\left\{1,k\right\}$.
Then $\frac{1}{\gcdd}\cdot D$ contains a multiple of $k$,
and yet $\chd\le3\le k$
(this is not hard to show; it will also follow from Lemma~\ref{L:dc3}).
As with general graphs, it appears to be quite difficult to
determine when a distance graph has a proper $k$-coloring, for $k\ge3$.
However, when $D$ is finite, there does exist an algorithm to
determine $\chd$.
This was proven for $D$ a finite set of primes
by Eggleton, Erd\H os, and Skilton
\cite[Corollary to Thm.~2]{EES90};
essentially the same proof works for more general sets.

\begin{theorem} \label{T:existsalg}
There exists an algorithm to determine $\chd$
for $D$ a finite set of positive integers.\end{theorem}

\begin{proof} (Outline---see~\cite[Thm.~2]{EES90})
Let $q=\max(D)$.
Then $\chd\le q+1$, by Lemma~\ref{L:nomult}.
We consider the colorings of the subgraph of $\gzd$
induced by $S=\iv{1}{q^q+q}$.
We show that, for $k\le q$, if $S$ has a proper
$k$-coloring, then $\chd\le k$;
thus, $\chd$ can be determined by a bounded search.

Let $k\le q$, and suppose that $S$ has a proper $k$-coloring.
The number of $k$-colorings of a block of $q$
consecutive integers is at most $q^q$.
Since $S$ contains $q^q+1$ such blocks,
two such blocks contained in $S$
(say $\iv{a}{a+q-1}$ and $\iv{b}{b+q-1}$, with $a<b$)
receive the same pattern of colors.
We extend the coloring of $\iv{a}{b+q-1}$ to a
coloring $f$ of $\bZ$ using the rule
$f(i+a-b)=f(i)$, for all $i$.
We can show that this is a proper coloring if $\gzd$,
and so $\chd\le k$.\ggcendpf\end{proof}

While an algorithm exists to determine $\chd$ for
finite $D$, we do not know whether
there is an efficient algorithm.
For finite graphs, determining whether the chromatic number is
at most $k$ is NP-complete~\cite{GaJo79}.
We conjecture that this is also true for distance graphs
with finite distance sets.

\begin{conjecture} \label{J:3npc}
Let $k\ge3$.
Determining whether $\chd\le k$ for finite sets $D$
is NP-complete.\ggcnopf\end{conjecture}

\section{Divisibility Chains} \label{S:dc}

We now focus on a particular class of distance graphs:
those in which the distance set is
totally ordered by divisibility.
We show that all such graphs are $4$-colorable,
and we determine their chromatic numbers.

A \emph{divisibility chain} is a set of positive integers that is
totally ordered by the divisibility relation.
When $D$ is a (finite or infinite) divisibility chain we
denote the elements of $D$ by $d_1,d_2,\dotsc$, where
$d_1\mid d_2\mid\dotsb$.
The \emph{ratios} of $D$ are the numbers
$r_i=\frac{d_{i+1}}{d_i}$, for each $i$.
When determining $\chd$, we may, by Lemma~\ref{L:conn},
assume that $\gcdd=d_1=1$.
Thus, $\chd$ depends only on the ratios.
We may also assume that all the $d_i$'s are distinct, that is,
that none of the ratios is equal to $1$.

A \emph{string} over $\{1,2\}$ is a finite sequence of $1$'s
and $2$'s, written without spaces or separators.
For example, $\alpha=1211$ is a string of length 4
with $\alpha_1=1$, $\alpha_2=2$, etc.

For $k$ a positive integer, a string $\alpha$
is \emph{$k$-compatible} with a distance set $D$
if there is a proper $k$-coloring of $\gzd$ with
colors $\iv{0}{k-1}$ such that the differences, modulo $k$,
between colors of consecutive vertices form
repeated copies of $\alpha$.
Below is part of such a coloring with $k=4$ and $\alpha=1211$.
\newcommand{\qs}{{}}
\newcommand{\qp}{\phantom{0}}
\newcommand{\qx}{{\,}}
\newcommand{\qpqpqs}{\qp\qp\qs}
\newcommand{\qxqpqs}{\qx\qp\qs}
\[
\begin{array}{rl}
\text{vertex}&
  \qx0\qs \qpqpqs \qp1\qs \qpqpqs \qp2\qs \qpqpqs \qp3\qs \qpqpqs
  \qp4\qs \qpqpqs \qp5\qs \qpqpqs \qp6\qs \qpqpqs \qp7\qs \qpqpqs
  \qp8\qs \qpqpqs \qp9\qs \qpqpqs 10\qs   \qpqpqs 11\qs   \qpqpqs
  12\qs   \cr
\text{color}&
  \qx0\qs \qpqpqs \qp1\qs \qpqpqs \qp3\qs \qpqpqs \qp0\qs \qpqpqs
  \qp1\qs \qpqpqs \qp2\qs \qpqpqs \qp0\qs \qpqpqs \qp1\qs \qpqpqs
  \qp2\qs \qpqpqs \qp3\qs \qpqpqs \qp1\qs \qpqpqs \qp2\qs \qpqpqs
  \qp3\qs \cr
\text{difference}&
  \qxqpqs \qp1\qs \qpqpqs \qp2\qs \qpqpqs \qp1\qs \qpqpqs \qp1\qs
  \qpqpqs \qp1\qs \qpqpqs \qp2\qs \qpqpqs \qp1\qs \qpqpqs \qp1\qs
  \qpqpqs \qp1\qs \qpqpqs \qp2\qs \qpqpqs \qp1\qs \qpqpqs \qp1\qs
  \qpqpqs \cr
\end{array}
\]
We see that $1211$ is not $4$-compatible with $\{3\}$,
since, for example, $1$ and $4$ receive the same color;
this is because the sum of three consecutive entries of
the repeated copies of $\alpha$ is divisible
by $4$ (i.e., $2+1+1=4$).
Generally,
a string $\alpha$ is $k$-compatible with $\{d\}$ if the
concatenation of repeated copies of $\alpha$ contains no
$d$ consecutive entries whose sum is a multiple of $k$.

\begin{theorem} \label{T:dc4}
If $D$ is a divisibility chain,
then $\chd\le4$.
\end{theorem}

\begin{proof}
We may assume that the $d_i$'s are all distinct, and that $d_1=1$.
We use notation such as $\alpha^n$ to denote a string;
the superscript does not denote exponentiation or concatenation.

\smallskip\noindent\emph{Claim.}
For $n=1,2,3,\dotsc$, there exist strings $\alpha^n$,
$\beta^n$ of length $d_n$ over $\left\{1,2\right\}$ such that
\begin{enumerate}
\item $\alpha^n$, $\beta^n$ differ only in the first entry,
with $\alpha^n_1=1$, and $\beta^n_1=2$, and
  \label{I:dc4pf-differ}
\item if $\gamma$ is a string resulting from the
concatenation of any number of copies of $\alpha^n$ and/or
$\beta^n$, in any order, then $\gamma$ is $4$-compatible with
$\left\{d_1,d_2,d_3,\dotsc,d_n\right\}$.
  \label{I:dc4pf-compat}
\end{enumerate}

Before we prove the claim, we show that the theorem follows from it.
If the claim holds, then, for each $n$,
$\alpha^n$ is $4$-compatible with
$\left\{d_1,d_2,\dotsc,d_n\right\}$, and so
$\gzmatch{\left\{d_1,d_2,\dotsc,d_n\right\}}$
has a proper 4-coloring.
Since every finite subgraph of $\gzd$
is isomorphic to a finite subgraph of
$\gzmatch{\left\{d_1,\dotsc,d_n\right\}}$ for some $n$,
every finite subgraph of $\gzd$ is 4-colorable,
and we may conclude that $\chd\le4$, by a compactness
argument.
Hence, it suffices to prove the claim.

\smallskip\noindent\emph{Proof of Claim.}
We proceed by induction on $n$.
For $n=1$, we assumed that $d_n=1$.
Let $\alpha^n=1$, and let $\beta^n=2$;
these satisfy the claim for $n=1$.

Now suppose that $n\ge1$,
and that the claim holds for $n$.
Define $s$ and $t$ as follows.
\[
s:=\sum_{i=1}^{d_n}\alpha^n_i;\qquad
t:=\sum_{i=1}^{d_n}\beta^n_i=s+1.
\]

We show first that $\iv{r_n\cdot s}{r_n\cdot t}$ contains
integers $w$, $w+1$, neither a multiple of $4$.
If $r_n>2$, then this is true since there are at
least 4 consecutive integers in $\iv{r_n\cdot s}{r_n\cdot t}$.
On the other hand, if $r_n=2$,
then $r_n\cdot s$ and $r_n\cdot t$ are both even.
Exactly one of the two is divisible by four.
If $4\mid\left(r_n\cdot s\right)$, then let $w=r_n\cdot s+1$;
otherwise, let $w=r_n\cdot s$.

Now we choose $a\ge1$, $b\ge0$ so that $a+b=r_n$
and $as+bt=w$:
let $b=w-r_n\cdot s$, and let $a=r_n-b$.
We define $\alpha^{n+1}$ to be the concatenation of
$a$ copies of
$\alpha^n$ followed by $b$ copies of $\beta^n$.
We let
$\beta^{n+1}$ be the concatenation of $\beta^n$
followed by $a-1$ copies of $\alpha^n$ followed by $b$
copies of $\beta^n$;
equivalently,
$\beta^{n+1}$ is $\alpha^{n+1}$ with its first
entry replaced by $2$.

Now, $\alpha^{n+1}$ and $\beta^{n+1}$ both have length $d_{n+1}$,
since $a+b=r_n$, and
$\alpha^{n+1}$ and $\beta^{n+1}$ differ only in the first entry.
Let $\gamma$ be a concatenation of copies of $\alpha^{n+1}$, $\beta^{n+1}$.
Then $\gamma$ is a concatenation of copies of $\alpha^n$ and $\beta^n$,
and so, by the induction hypothesis, $\gamma$ is $4$-compatible
with $\left\{d_1,d_2,\dotsc,d_n\right\}$.

In order to prove that $\alpha^{n+1}$, $\beta^{n+1}$ satisfy the
claim, it remains only to show that $\gamma$ is $4$-compatible
with $\left\{d_{n+1}\right\}$.
This is true if the concatenation of repeated copies of $\gamma$
has no $d_{n+1}$ consecutive entries whose sum is a multiple
of $4$.
Since $\alpha^{n+1}$ and $\beta^{n+1}$ differ in only one entry,
the sum of $d_{n+1}$ consecutive entries of
repeated copies of $\gamma$ is equal either to the sum
of the entries of $\alpha^{n+1}$ or to the sum of the entries
of $\beta^{n+1}$;
that is, it is equal
\begin{align*}
\text{either to}\quad\sum_{i=1}^{d_{n+1}}\alpha^{n+1}_i&=as+bt=w,\\
\text{or to}\quad\sum_{i=1}^{d_{n+1}}\beta^{n+1}_i&=t+(a-1)s+bt=w+1.
\end{align*}
Neither of these is a multiple of $4$.

Thus, the claim is proven.\ggcendpf\end{proof}

The bound in Theorem~\ref{T:dc4} is sharp:
if $D=\left\{1,2,6\right\}$,
then the subgraph of $\gzd$ induced by
$\iv{1}{7}$ has no proper 3-coloring.
On the other hand, graphs satisfying the hypotheses of
the theorem need not have chromatic number 4, even if
$D$ is infinite.
For example, if $D$ is the set of all powers of $3$, then
every element of $D$ is odd, and so $\chd=2$,
by Proposition~\ref{P:chi2}.

\begin{example} \label{G:usedc4}
Let $D=\left\{d_1,d_2,d_3,\dotsc\right\}$, where
$d_i=i!$ for each $i$.
We use the technique of the above proof to produce part of a
proper 4-coloring of $\gzd$, the factorial distance graph.

Let $\alpha^1=1$ and $\beta^1=2$.
We find consecutive nonmultiples of 4
in $\iv{2\cdot1}{2\cdot2}=\{2,3,4\}$:
let $w=2$, so that $w+1=3$.
So, $a=2$, and $b=0$.
The string $\alpha^2$ is $2$ copies of $\alpha^1$
followed by $0$ copies of $\beta^1$.
That is, $\alpha^2=11$, and so $\beta^2=21$.

Continuing, we find consecutive nonmultiples of 4
in $\iv{3\cdot2}{3\cdot3}=\{6,7,8,9\}$:
let $w=6$, so that $w+1=7$.
So, $a=3$, and $b=0$.
The string $\alpha^3$ is $3$ copies of $\alpha^2$ followed by
$0$ copies of $\beta^2$.
That is, $\alpha^3=111111$, and so $\beta ^3=211111$.

Once again, we find consecutive nonmultiples of 4
in $\iv{4\cdot6}{4\cdot7}=\iv{24}{28}$:
let $w=25$, so that $w+1=26$.
So, $a=3$, and $b=1$.
The string $\alpha^4$ is $3$ copies of $\alpha^3$ followed
by $1$ copy of $\beta^3$.
That is, $\alpha^4=111111111111111111211111$,
and so $\beta^4=211111111111111111211111$.

The coloring of $\iv{1}{24}$ obtained from $\alpha^4$ is the following.
\[0,1,2,3,0,1,2,3,0,1,2,3,0,1,2,3,0,1,2,0,1,2,3,0.\ggcendeg\]
\end{example}

In \emph{almost} the entire proof of Theorem~\ref{T:dc4},
``4'' can be replaced by ``3'';
that is, we use $3$-compatibility instead of
$4$-compatibility,
we find a $3$-coloring instead of a $4$-coloring,
and we find consecutive nonmultiples of $3$ instead of $4$.
The one place where $4$ is required is the
argument in the proof showing the existence of
two consecutive nonmultiples when $r_n=2$.
Thus, if we require that $r_n\ne 2$ for each $n$,
then we can replace $4$ by $3$ in the proof, and
we have the following result.

\begin{lemma} \label{L:dc3}
Let $D$ be a divisibility chain, with ratios $r_1,r_2,\dotsc$.
If $r_i\ne2$ for all $i$,
then $\chd\le3$.\ggcnopf\end{lemma}

Again, the bound in this result is sharp:
if $D=\left\{1,4\right\}$, then the subgraph of $\gzd$
induced by $\iv{1}{5}$ has no proper 2-coloring.

We now find $\chd$ for every divisibility chain
$D$.

\begin{theorem} \label{T:dcchi}
Let $D$ be a divisibility chain, with ratios $r_1,r_2,\dotsc$.
All of the following hold.
\begin{enumerate}
\item $\chd\le4$.
  \label{I:dcchi4}
\item $\chd\le3$ if and only if there do not exist
  $i$, $j$ with $i<j$, $r_i=2$, and $3\mid r_j$.
  \label{I:dcchi3}
\item $\chd\le2$ if and only if $r_i$ is odd,
  for each $i$.
  \label{I:dcchi2}
\item $\chd=1$ if and only if $D=\emptyset$.
  \label{I:dcchi1}
\end{enumerate}
\end{theorem}

\begin{proof}
Statement~\ref{I:dcchi4} follows from Theorem~\ref{T:dc4},
statement~\ref{I:dcchi2} follows from Proposition~\ref{P:chi2},
and statement~\ref{I:dcchi1} holds because a graph is 1-colorable
precisely when it has no edges.
It remains to prove statement~\ref{I:dcchi3}.
We may assume that $d_1=1$.

$(\Longrightarrow)$
Suppose that there exist $i$ and $j$ with $i<j$, $r_i=2$, and $3\mid r_j$.
Then $d_{i+1}=2d_i$, and $d_{j+1}$ is divisible by $3d_i$.
Suppose that $\gzd$ has a proper $3$-coloring.
Consider the colors assigned to the multiples of $d_i$.
Since $d_i,2d_i\in D$, vertices $0$, $d_i$, and $2d_i$
induce a complete subgraph and so must be assigned 3 different
colors.
Similarly, $d_i$, $2d_i$, and $3d_i$ must receive 3 different
colors, and so $0$ and $3d_i$ have the same color.
Continuing this argument, all multiples of $3d_i$ must receive the
same color, including $0$ and $d_{j+1}$, which is impossible.

$(\Longleftarrow)$
Suppose there do not exist $i$ and $j$ with the properties
specified in statement~\ref{I:dcchi3};
that is, every ratio divisible by $3$ precedes every ratio equal to
$2$ in the list $\left\{r_1,r_2,\dotsc\right\}$.
If there exist infinitely many ratios that are divisible by $3$,
then, by our assumption, there exists no ratio equal to $2$,
and so $\chd\le3$, by Lemma~\ref{L:dc3}.
Thus, we may assume that there are only finitely many
ratios that are divisible by $3$.

Let $c$ be the least positive integer such that
$3\nmid r_i$, for all $i\ge c$.
Then none of $r_1,r_2,\dotsc,r_{c-1}$ is equal to $2$.
Thus, by Lemma~\ref{L:dc3}, the graph
$\gzmatch{\left\{d_1,d_2,\dotsc,d_c\right\}}$
has a proper $3$-coloring.
By the proof of Lemma~\ref{L:dc3}---that is, the proof of
Theorem~\ref{T:dc4}, as modified to prove Lemma~\ref{L:dc3}---there
is a string $\alpha^c$ of length $d_c$
over $\left\{1,2\right\}$ such that
$\alpha^c$ is $3$-compatible with
$\left\{d_1,d_2,\dotsc,d_c\right\}$.

We claim that $\alpha^c$ is $3$-compatible with $D$.
To see this, first note that
\[
\sum_{i=1}^{d_c}\alpha^c_i
\]
is not a multiple of $3$, since $\alpha^c$ is $3$-compatible
with $\left\{d_c\right\}$.
Thus, if integers $x$ and $y$ differ by a multiple of $d_c$,
then, in a $3$-coloring whose differences, modulo $3$,
form repeated copies of $\alpha^c$,
$x$ and $y$ receive the same color
precisely when their difference is a multiple of $3d_c$.
Now, no $r_i$ with $i\ge c$ is divisible by $3$;
thus, no $d_i$ with $i\ge c$ is divisible by $3d_c$.
We conclude that, for each $i\ge c$,
no two integers with difference $d_i$ receive the same color,
and so $\alpha^c$ is $3$-compatible with
$\left\{d_c,d_{c+1},d_{c+2},\dotsc\right\}$.

Thus, $\alpha^c$ is $3$-compatible with $D$, and we have
$\chd\le3$.\ggcendpf\end{proof}

By Theorem~\ref{T:dcchi}, the chromatic number of
the factorial distance graph is $4$.
We will have more to say about this graph in the next section.

\section{Periodic Colorings} \label{S:pc}

In this section, we consider periodic proper colorings
of distance graphs.
We characterize those distance graphs that have no periodic
proper coloring,
and we find a relationship between
the chromatic number and the nonexistence of periodic proper
colorings.
Periodic colorings have been previously studied in~\cite{EES90}.

\begin{lemma} \label{L:nomultper}
Let $D$ be a set of positive integers, and
let $k$ be a positive integer.
If $D$ contains no multiple of $k$,
then $\gzd$
has a periodic proper $k$-coloring.
\end{lemma}

\begin{proof}
We may assume $D\ne\emptyset$.
The proof of Lemma~\ref{L:nomult}
gives a periodic proper $k$-coloring of each component of
$\gzd$;
this results in a periodic proper $k$-coloring
of the graph.\ggcendpf\end{proof}

We can use Lemma~\ref{L:nomultper} to characterize
those distance graphs that have no periodic proper coloring.
The following result generalizes an observation of
Eggleton~\cite{EggR97}
that the square distance graph has no periodic proper coloring.

\begin{proposition} \label{P:noper}
Let $D$ be a set of positive integers.
The graph $\gzd$ has no periodic proper coloring
if and only if $D$ contains a multiple of every positive integer.
\end{proposition}

\begin{proof}
$(\Longrightarrow)$
If there is some positive integer $k$ such
that $D$ contains no multiple of $k$, then,
by Lemma~\ref{L:nomultper},
$\gzd$ has a periodic proper coloring.

$(\Longleftarrow)$
Let $D$ contain a multiple of every positive integer.
Let $\gzd$ be colored in a periodic manner;
say this coloring has period $k$.
Every pair of vertices whose difference is a multiple
of $k$ will have the same color.
Since $D$ contains some multiple of $k$, this cannot be a
proper coloring.\ggcendpf\end{proof}

\begin{remark} \label{R:fact}
It follows from Theorem~\ref{T:dcchi}
that the chromatic number of the factorial distance graph
is $4$.
However, by Proposition~\ref{P:noper}, the factorial
distance graph has no periodic proper coloring.\ggcendrem\end{remark}

Now we examine the effect of the existence of uniquely colorable
subgraphs on proper colorings of distance graphs.
We prove a useful lower bound on the chromatic number
based on uniquely colorable subgraphs and periodic
colorings.

\begin{proposition} \label{P:uniqueper}
Let $D$ be a set of positive integers, and
let $k$ be a positive integer.
If $\gzd$ has a finite,
uniquely $k$-colorable subgraph, then
every proper $k$-coloring of $\gzd$ is periodic.
\end{proposition}

\begin{proof}
Suppose that $H$ is a uniquely $k$-colorable subgraph of
$\gzd$.
We may assume that the least integer that is a vertex of $H$ is $1$.
Let $n$ be the greatest-numbered vertex of $H$.
Since $H$ is uniquely $k$-colorable, every $k$-coloring of
$\iv{1}{n-1}$ that can be extended to a proper $k$-coloring of
$\gzd$ has a unique extension to a proper $k$-coloring
of $\iv{1}{n}$.

In short, once we have $k$-colored $\iv{1}{n-1}$, the
color of vertex $n$ is forced.
But, $\iv{2}{n+1}$ also contains a copy of $H$, and so
once we have colored $\iv{2}{n}$, the color of vertex $n+1$
is forced.
By an inductive argument, we can see that $k$-coloring
$\iv{1}{n-1}$ completely determines the coloring of
$\left[1,\infty\right)$.

Essentially the same argument works in the opposite direction:
$k$-coloring $\iv{2}{n}$ forces a certain color to occur
at vertex $1$.
Hence, $k$-coloring any set of $n-1$ consecutive vertices determines
the coloring of all of $\bZ$.

Now, there are only a finite number of $k$-colorings of $n-1$
consecutive integers.
Since the colorings of blocks of $n-1$ consecutive integers
must eventually repeat, every proper $k$-coloring of the distance
graph is periodic.\ggcendpf\end{proof}

Proposition~\ref{P:noper} and Proposition~\ref{P:uniqueper}
have nearly opposite conclusions;
the former concludes that the graph has
no periodic proper coloring,
while
the latter concludes that every proper $k$-coloring of the
distance graph is periodic.
Suppose that a distance graph satisfies the hypothesis of
both propositions, that is,
the distance set contains a multiple of every
positive integer,
and
the graph has a finite, uniquely $k$-colorable subgraph.
Then the conclusions of both propositions must be true:
there is no periodic proper coloring,
and yet
every proper $k$-coloring is periodic.
We can only conclude that the distance graph must have no
proper $k$-coloring at all, and so we have the following result.

\begin{theorem} \label{T:omega+1}
Let $D$ be a set of positive integers, and
let $k$ be a positive integer.
If $D$ contains a multiple of every positive integer,
and $\gzd$ has a finite, uniquely $k$-colorable subgraph,
then $\chd\ge k+1$.\ggcnopf\end{theorem}

We can use Theorem~\ref{T:omega+1} to place a lower bound on
the chromatic number of the square distance graph.
Let $D$ be the set of all positive squares.
Any Pythagorean triple gives a $K_3$ in the square distance graph.
For example, the vertices $0$, $3^2$, $5^2$ induce a $K_3$, since
$3,4,5$ is a Pythagorean triple.
Since $\gzd$ has a $K_3$ subgraph,
$\chd\ge3$.
Furthermore, $K_3$ is uniquely $3$-colorable, and $D$ contains
a multiple of every positive integer.
Thus, $\chd\ge4$, by Theorem~\ref{T:omega+1}.
Eggleton~\cite{EggR97} has found a $K_4$ in the square distance graph:
the vertices are $0$, $672^2$, $680^2$, and $697^2$.
We have $680^2-672^2=104^2$, $697^2-680^2=153^2$, and
$697^2-672^2=185^2$.
Since the square distance graph has a uniquely $4$-colorable
subgraph, we have the following result.

\begin{corollary} \label{C:sqchige5}
The chromatic number of the square distance graph is
at least $5$.\ggcnopf\end{corollary}

We do not know whether
the square distance graph contains a $K_5$ or whether its
chromatic number is greater than $5$.

\begin{problem} \label{O:sqchi}
What is the chromatic number of the square distance graph?
Equivalently, what is the least number of classes into
which the integers can be partitioned,
so that no two members of the same class differ by a
square?\ggcnopf\end{problem}

It seems likely that no finite number of colors suffices.

We can ask similar questions about the distance graph
resulting when $D$ is the set of all positive $n$th
powers, for $n\ge3$.
We know that these graphs contain no $K_3$
(this is equivalent to ``Fermat's Last Theorem'',
proven by Wiles \cite{WilA95}),
that they do not have periodic proper colorings,
by Proposition~\ref{P:noper},
and that their chromatic numbers are all at least $3$, by
Theorem~\ref{T:omega+1} (or Proposition~\ref{P:chi2}).
It seems likely that these graphs have infinite chromatic
number as well.

As noted in Section~\ref{S:intro},
determining which distance graphs have chromatic number
at most $k$, for a given $k\ge3$, appears to be
difficult.
A similar problem, whose difficultly we cannot estimate
at this time, is the following.

\begin{problem} \label{O:charinf}
Characterize those sets $D$ such that
$\chd$ is\break
infinite.\ggcnopf
\end{problem}

No coloring requiring an infinite number of colors is periodic.
Thus, by Proposition~\ref{P:noper}, a necessary condition for
such sets $D$ is that they contain a multiple of every integer.
However, this condition is not sufficient, by Remark~\ref{R:fact}.

\section*{Acknowledgments} \label{S:ack}

The author is grateful to Professor Roger Eggleton for bringing
this topic to his attention and for helpful discussions.

\enddocument